\documentclass[sn-vancouver,Numbered]{sn-jnl}

\usepackage{graphicx}%
\usepackage{multirow}%
\usepackage{amsmath,amssymb,amsfonts}%
\usepackage{amsthm}%
\usepackage{mathrsfs}%
\usepackage[title]{appendix}%
\usepackage{xcolor}%
\usepackage{textcomp}%
\usepackage{manyfoot}%
\usepackage{booktabs}%
\usepackage{algorithm}%
\usepackage{algorithmicx}%
\usepackage{algpseudocode}%
\usepackage{listings}%
\usepackage{anyfontsize} 
\newcommand{\R}{\mathbb R}           
\newcommand{\bfX}{{\bf X}}
\usepackage{lineno}
\usepackage[version=4]{mhchem}

\newcommand{\insitu}{{\emph{in situ}}}
\newcommand{\bigO}{\mathcal{O}}

\begin{document}

The submitted manuscript has been created by UChicago Argonne, LLC, Operator of Argonne National Laboratory(“Argonne”). Argonne, a U.S. Department of Energy Office of Science laboratory, is operated under Contract No. DE-AC02-06CH11357. The U.S. Government retains for itself, and others acting on its behalf, a paid-up nonexclusive, irrevocable worldwide license in said article to reproduce, prepare derivative works, distribute copies to the public, and perform publicly and display publicly, by or on behalf of the Government.
\newpage

\title[Article Title]{Homomorphic data compression for real time photon correlation analysis}


\author[1]{\fnm{Sebastian} \sur{Strempfer}}\nomail
\author[2]{\fnm{Zichao Wendy} \sur{Di}}\nomail
\author[2]{\fnm{Kazutomo} \sur{Yoshii}}\nomail
\author[3]{\fnm{Yue} \sur{Cao}}\nomail
\author[1]{\fnm{Qingteng} \sur{Zhang}}\nomail
\author[1]{\fnm{Eric M.} \sur{Dufresne}}\nomail
\author[1]{\fnm{Mathew} \sur{Cherukara}}\nomail
\author[1]{\fnm{Suresh} \sur{Narayanan}}\nomail
\author[4]{\fnm{Martin V.} \sur{Holt}}\nomail
\author*[1]{\fnm{Antonino} \sur{Miceli}}\email{amiceli@anl.gov}
\author*[4]{\fnm{Tao} \sur{Zhou}}\email{tzhou@anl.gov}

\affil[1]{\orgdiv{Advanced Photon Source}, \orgname{Argonne National Laboratory}, \orgaddress{\street{9700 S Cass Ave}, \city{Lemont}, \postcode{60439}, \state{IL}, \country{USA}}}

\affil[2]{\orgdiv{Mathematics and Computer Science}, \orgname{Argonne National Laboratory}, \orgaddress{\street{9700 S Cass Ave}, \city{Lemont}, \postcode{60439}, \state{IL}, \country{USA}}}

\affil[3]{\orgdiv{Materials Science Division}, \orgname{Argonne National Laboratory}, \orgaddress{\street{9700 S Cass Ave}, \city{Lemont}, \postcode{60439}, \state{IL}, \country{USA}}}

\affil[4]{\orgdiv{Center for Nanoscale Materials}, \orgname{Argonne National Laboratory}, \orgaddress{\street{9700 S Cass Ave}, \city{Lemont}, \postcode{60439}, \state{IL}, \country{USA}}}


\abstract{The construction of highly coherent x-ray sources has enabled new research opportunities across the scientific landscape. The maximum raw data rate per beamline now exceeds 40 GB/s, posing unprecedented challenges for the online processing and offline storage of the big data. Such challenge is particularly prominent for x-ray photon correlation spectroscopy (XPCS), where real time analyses require simultaneous calculation on all the previously acquired data in the time series. We present a homomorphic compression scheme to effectively reduce the computational time and memory space required for XPCS analysis. Leveraging similarities in the mathematical expression between a matrix-based compression algorithm and the correlation calculation, our approach allows direct operation on the compressed data without their decompression. The lossy compression reduces the computational time by a factor of 10,000, enabling real time calculation of the correlation functions at kHz framerate.
Our demonstration of a homomorphic compression of scientific data provides an effective solution to the big data challenge at coherent light sources. Beyond the example shown in this work, 
the framework can be extended to facilitate real-time operations directly on a compressed data stream for other techniques.}




\maketitle

\section{Introduction}\label{sec1}

The construction of highly coherent x-ray sources around the world has enabled new research opportunities across the scientific landscape\cite{Raimondi2023, Martensson2018, APS-U}. The 100-fold increase in source brightness also brings about a drastic increase in the raw data rate. The maximum raw data rate per beamline exceeds 40 GB/s or 20 PB/week, creating unprecedented challenges for both online processing and offline storage of the big data\cite{Eriksson2014}. Among those that benefit the most from the new sources are a number of techniques that exploits the high coherence, namely coherent diffraction imaging\cite{Clark2013, Ulvestad2015}, ptychography\cite{Thibault2008, Hruszkewycz2017}, and x-ray photon correlation spectroscopy (XPCS)\cite{Shpyrko2014, Sinha2014}. XPCS is a powerful technique for studying the dynamics of materials at the nanoscale\cite{Cao2020}. In a typical XPCS experiment, a megapixel two-dimensional detector is used to record the speckle pattern from a coherent x-ray beam scattered from the sample. By measuring the intensity fluctuations of the speckle pattern as a function of time, the dynamics of various equilibrium and non-equilibrium processes can be studied with a picoseconds-to-seconds temporal resolution\cite{Zhang2018, Seaberg2017, Roseker2018}. Unlike most other techniques where real time analysis is performed only on the most current image (data)\cite{Liu2021, Babu2023}, the real time analysis in XPCS is highly compute-intensive as it requires simultaneous calculation on all the previously acquired images in the time series. Previous demonstration of high-speed XPCS shows a raw data rate of 6.5 GB/s on a 500k pixelated detector running at 52 kHz\cite{Zhang2021}. The maximum raw data rate is thus 40 GB/s with the latest generation of 3 megapixel detector\cite{Rigaku}. It would require tens of GB of disk space to store the raw data generated in a mere second of measurement, and hundreds of TB of computer memory for the vectorized calculation of the time-correlation function that exceeds the capacity of most local computing. More robust computational capacity is available at high performance computing (HPC) facilities\cite{Salim2019}, although the lengthy data transfer process and competition for on-demand resources prohibit real-time analysis and feedback for experimental steering\cite{Kandel2023}.

Big data generated by an XPCS experiment is typically compressed\cite{areadetector} or sparsified\cite{Zhang2021} before being written to the disk. While state-of-the-art lossless\cite{Welch1984} and lossy\cite{Zhao2021} compression algorithms can effectively reduce the size of the big data to be stored offline, the decompression of the big data during the analysis stage still poses a huge challenge for both the computational resource and the analysis program\cite{Chu2022}. In this work, we present a homomorphic compression scheme that enables direct operation on compressed XPCS data without their decompression. The approach is made possible by leveraging similarities in the mathematical expression between singular value decomposition (SVD) matrix-based compression\cite{Mazurov2017} and the vectorized computation of the two-time correlation (TTC) function in XPCS\cite{Bikondoa2017}. With lossless compression, we show that the photon correlation of the compressed data is an exact match for the photon correlation of the original data, despite computations being 800 times faster and on a dataset 800 times smaller. With lossy compression, we achieve a compression ratio of 10,000 while preserving a qualitative approximation of the result using the original data. The lossy compression reduces the computational time by a factor of 10,000, which paves the way for in-pixel compression on the x-ray detector \cite{valentin, Strempfer2022} and real time calculation of the correlation functions at kHz framerate.

\section{Methods}\label{sec2}
\subsection{Direct Computation of TTC using SVD Compressed Data}\label{subsec21}

In an XPCS experiment, the two-time correlation function\cite{Bikondoa2017} is commonly defined as
\begin{equation}
\textbf{G}_{t_1,t_2} = \frac{\overline{I_{t_1}I_{t_2}}}{\overline{I_{t_1}}\;\overline{I_{t_2}}},
\end{equation}
where $I_{t_i}$ is the x-ray intensity measured at time $t_i$. $\textbf{G}\in \R^{N \times N}$ with $N$ being the number of images in the time series. In terms of matrix operation, $\textbf{G}$ can be computed as 
\begin{equation}
\textbf{G} = \textbf{X}\textbf{X}^\text{T},
\label{eq:2-tcf}
\end{equation}
where $\textbf{X}\in \R^{N \times M}$ is the XPCS data of the entire time series (Fig.~\ref{fig:workflow}a). $M$ is the number of total pixels per image. The superscript $^\text{T}$ denotes the transpose operator. Row $i$ of the matrix $\textbf{X}$ is an array of size $1 \times M$, corresponding to the flattened array of the $i^\text{th}$ image of the time series $I_{t_i}$ normalized by its Frobenius norm. 

For the singular value decomposition of data \textbf{X}, 
\begin{equation}
\textbf{X} = \textbf{U} \boldsymbol{\Sigma} \textbf{V}^\text{T}.
\label{eq:svd0}
\end{equation}
The columns of $\textbf{V}\in \R^{M \times N}$ are the eigenvectors of $\textbf{X}^\text{T}\textbf{X}$ and the diagonal entries of $\boldsymbol{\Sigma}\in \R^{N \times N}$ are the corresponding eigenvalues. The columns of $\textbf{U}\in \R^{N \times N}$ are the eigenvectors of $\textbf{X}\textbf{X}^\text{T}$. For lossless compressions, we multiply the data matrix by the entirety of \textbf{V} 
\begin{equation}
\textbf{Y} = \textbf{X} \textbf{V},
\label{eq:compress1}
\end{equation}
which encodes the original data $\textbf{X}\in \R^{N \times M}$ to its compressed form $\textbf{Y}\in \R^{N \times N}$; and
\begin{equation}
\bfX = \textbf{Y} \textbf{V}^\text{T} = \textbf{X} \textbf{V} \textbf{V}^\text{T}.
\label{eq:svd1}
\end{equation}
decodes \textbf{Y} back to the original data matrix \bfX. 

In the case of lossy compressions, the first $K$ eigenvectors of \textbf{V} corresponding to the highest $K$ eigenvalues in $\boldsymbol{\Sigma}$ are kept for the encoding of the data matrix
\begin{equation}
\textbf{V}_\text{K} \in \R^{M \times K} = \textbf{V}[:,1:K].
\label{eq:identity1}
\end{equation}
Then we have the lossy compressed data $\textbf{Y}_\text{K} \in \R^{N \times K}$ as
\begin{equation}
\textbf{Y}_\text{K} = \textbf{X} \textbf{V}_\text{K},
\label{eq:compress2}
\end{equation}
and the corresponding decompressed data $\tilde{\bfX} \in \R^{N \times M}$ as
\begin{equation}
\tilde{\bfX} = \textbf{X} \textbf{V}_\text{K} \textbf{V}_\text{K}^\text{T}.
\label{eq:svd2}
\end{equation}
$\textbf{X} = \tilde{\bfX}$ only if $K = \text{rank(\textbf{X})}$. We proceed to calculate an approximation $\tilde{\textbf{G}}\in \R^{N \times N}$ of the correlation \textbf{G} using the decompressed data matrix $\tilde{\bfX}$ and the original data matrix $\textbf{X}$. From Eq.~(\ref{eq:svd2}) we have
\begin{equation}
\tilde{\textbf{G}} = \tilde{\bfX} \textbf{X}^\text{T} = \textbf{X} \textbf{V}_\text{K}  \textbf{V}_\text{K}^\text{T} \textbf{X}^\text{T}.
\label{eq:final0}
\end{equation}
With Eqs.~(\ref{eq:compress2}) and (\ref{eq:final0}), it became obvious that for lossy compressions the approximated two-time correlation $\tilde{\textbf{G}}$ is none other than the two-time correlation function of the compressed matrix $\textbf{Y}_\text{K}$
\begin{equation}
\tilde{\textbf{G}} = \textbf{Y}_\text{K} \textbf{Y}_\text{K}^\text{T}.
\label{eq:final1}
\end{equation}
In other words, it is possible to compute an approximation of the two-time correlation of the raw data using the two-time correlation of SVD compressed data, without any decompression.

For lossless compression, because $\textbf{X} = \tilde{\textbf{X}}$, the two-time correlation of SVD compressed data is exactly the two-time correlation of the raw data
\begin{equation}
\textbf{G} = \textbf{X} \textbf{X}^\text{T} = \textbf{Y} \textbf{Y}^\text{T}= \tilde{\textbf{G}}.
\label{eq:final2}
\end{equation}

\subsection{Strategies for Offline and Online Compression}\label{subsec22}
\begin{figure}[htbp]
\centering\includegraphics[width=1\textwidth]{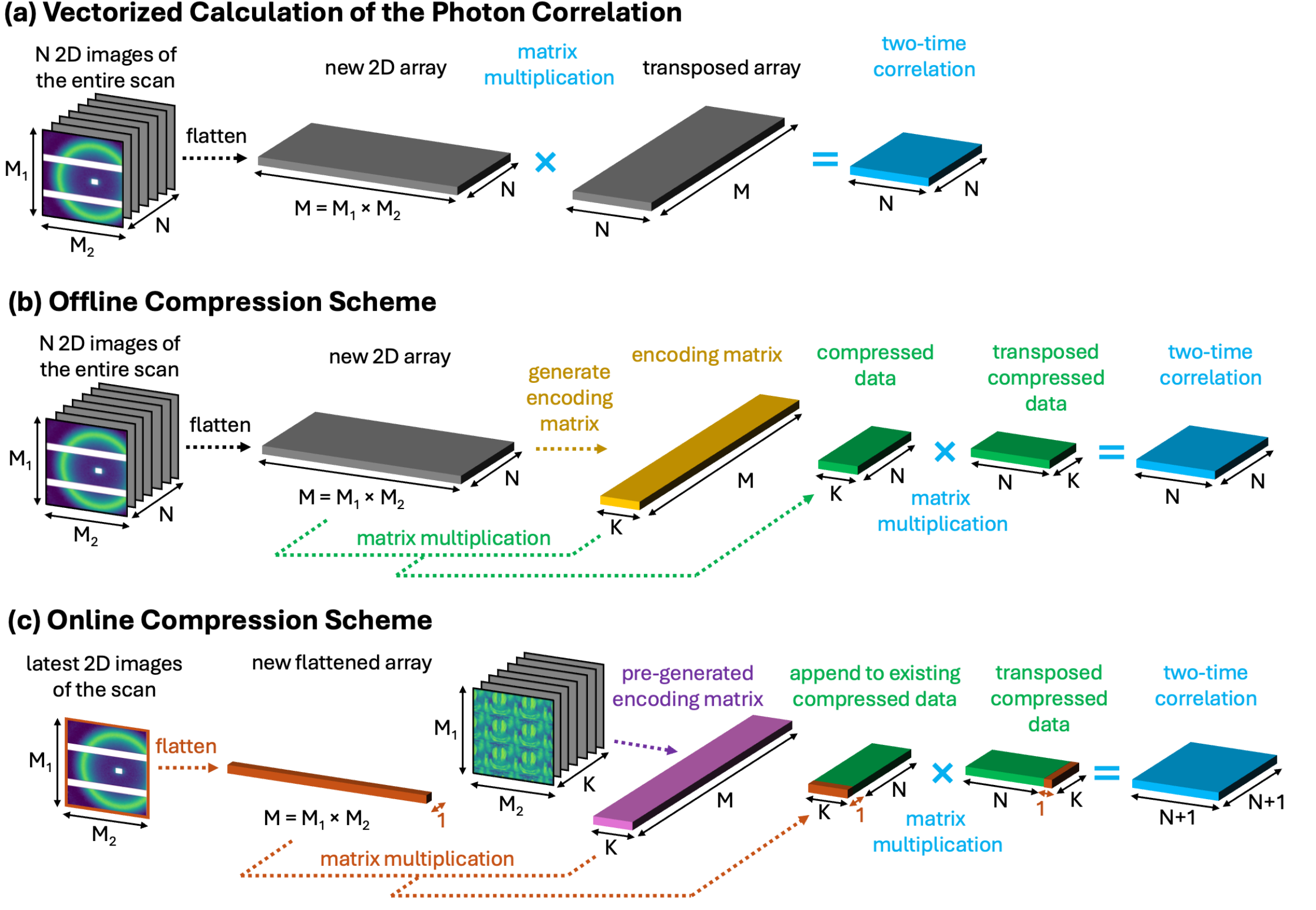}
\caption{Framework for XPCS analysis with (a) raw, (b) offline and (c) online compressed data.}
\label{fig:workflow}
\end{figure}

In Fig.~\ref{fig:workflow}b and c we present two strategies for the homomorphic compression and subsequent analysis of the XPCS data. For offline compression, the compression is initiated only at the end of the time series. The encoding matrix is generated by the data to-be-compressed itself, after flattening the two dimensions given by the detector frame. The data from the entire time series is compressed all at once, and the two-time correlation is computed solely with the compressed data. For online compression, the compression is performed on each image immediately after their individual acquisition. The encoding matrix is generated before the start of the time series, using data that may or may not be related to the experiment. In the illustrated example, we generate the encoding matrix using randomly shifted Mandrill image from the Signal and Image Processing Institute (SIPI) image database\cite{SIPI}. The encoding matrix generated in this manner has the advantage of being reusable in any experiment conducted at any facility. The newly compressed data takes the form of an 1D array, which is promptly appended to the previously compressed data of the time series. Because the size of the compressed data is significantly smaller than that of the raw data, this enables real-time computation of the two-time correlation of the XPCS time series. We note that it is also possible to compress and analyze a subset of the detector frame corresponding to pixels belonging to a specific $q$ range, in which case the data to be compressed is automatically a flattened array ready to be multiplied with the encoding matrix.

\section{Results}\label{sec3}
\subsection{High Fidelity with Offline Compression}\label{subsec31}

\begin{figure}[htbp]
\centering\includegraphics[width=1\textwidth]{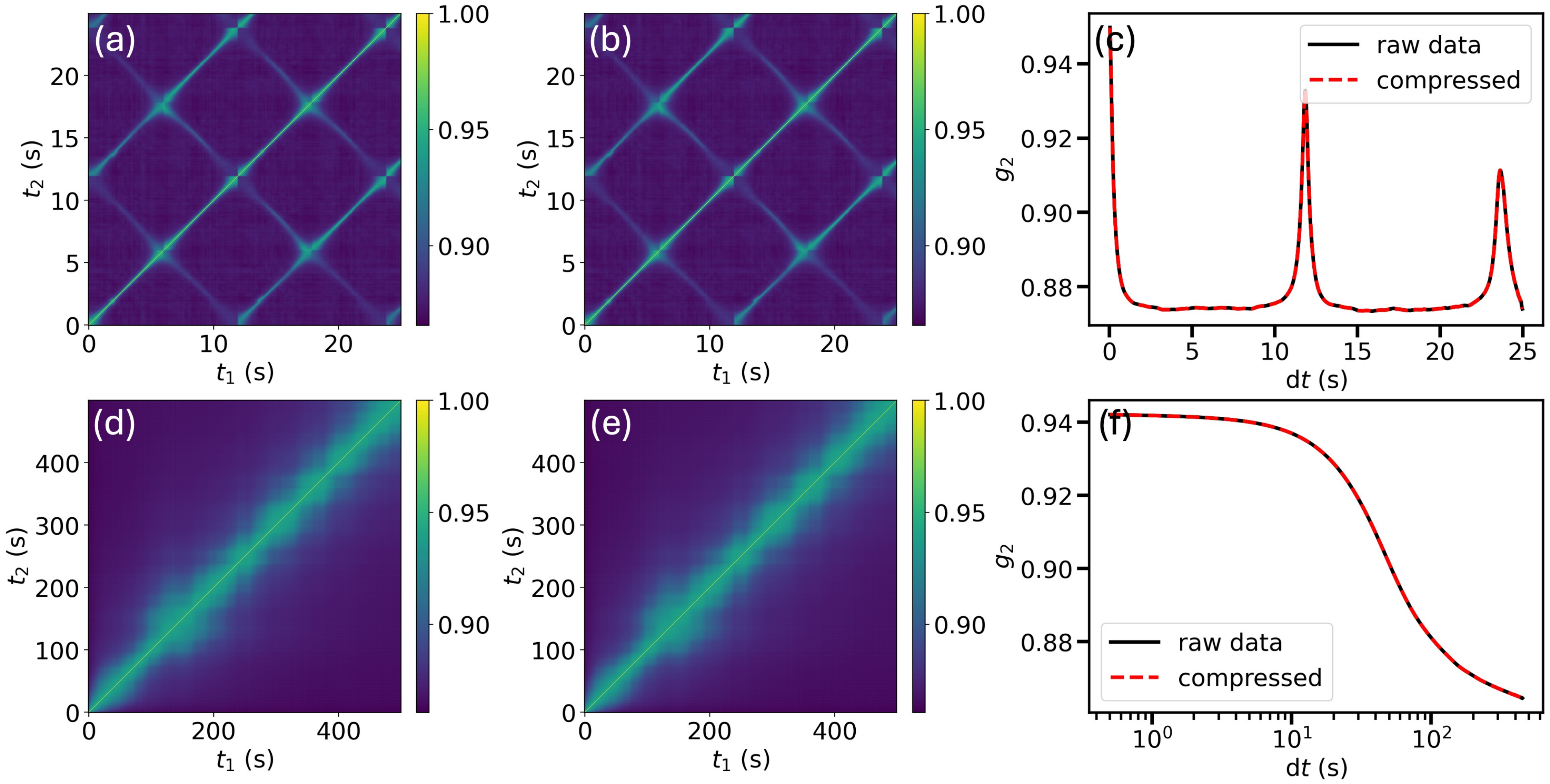}
\caption{Correlation calculation with offline lossless compression. (a) TTC of the raw oscillatory data. (b) TTC of the corresponding compressed data. (c) Comparison of $g_2$ between the raw and compressed data. (d) TTC of the raw rheology data. $M=610708$ pixels are extracted per image corresponding to the $q$ range between 0.002 and 0.003 $\textup{~\AA}^{-1}$. (e) TTC of the corresponding compressed data. (f) Comparison of $g_2$ between the raw and compressed data. }
\label{fig:offline_lossless}
\end{figure}
In the case of offline compression, we shall first explore the scenario of lossless compression. A significant advantage in this case, as shown in Eq.~(\ref{eq:final2}), is the possibility to compute the exact value of the raw data TTC using directly the compressed data. We demonstrate this using an experimental XPCS data on a silica filled rubber subjected to dynamic strain (hereinafter referred to as the oscillatory data)\cite{Presto2023}. The number of images in the time series is $N = 1000$ and the total number of pixel on the detector is $M=802896$. Fig.~\ref{fig:offline_lossless}a shows the TTC of the raw data, which equals the TTC of the compressed data shown in Fig.~\ref{fig:offline_lossless}b, despite the latter occupying a disk space that is 800 times smaller and executes 800 times faster. Fig.~\ref{fig:offline_lossless}c shows the value of the $g_2$ function calculated as 
\begin{equation}
g_2(q,\text{d}t) = \frac
{\langle I(q,t)\;I(q,t+\text{d}t)\rangle}
{\langle I(q) \rangle^2}.
\label{eq:g2}
\end{equation}
In this case, the entire $q$ range is used. The $g_2$ of the compressed data again equals the $g_2$ of the raw data.
In Fig.~\ref{fig:offline_lossless}d and Fig.~\ref{fig:offline_lossless}e we show respectively the TTC from another set of XPCS data on the \insitu{} stress relaxation of a silica colloidal glass (hereinafter referred to as the rheology data). In this case the raw data \textbf{X} is a subset ($M=610708$) of the detector pixels falling within a specific $q$ range. The number of images is still $N=1000$. Because SVD compression operates on a flattened 1D-array, it is directly applicable to the aforementioned subset despite it not being a regular 2D array. This is different from any transform coding based compression algorithm ({e.g.}, JPEG) which only operates on 2D images. The compressed data is capable to again reproducing a $g_2$ function that equals the $g_2$ of the raw data, as shown in Fig.~\ref{fig:offline_lossless}f.

We note that in the offline compression with $K = N$ ({i.e.}, lossless), we do not actually reduce the size of the data saved on the disk. This is because while the compressed data $\textbf{Y} \in \R^{N \times N}$ is about a few hundred times smaller than the raw data $\textbf{X} \in \R^{N \times M}$, one still needs to store the encoding matrix $\textbf{V} \in \R^{M \times N}$ which has the same size as the raw data. The effective total computation time is not always shorter either, due to the time necessary to generate $\textbf{V}$. Our proposed approach becomes advantageous when recalculating the TTCs, which frequently happens during subsequent data analysis. Once the compressed data are generated, they can be directly used for any future calculations. For each of those subsequent recalculations, the memory space requirement is tens of thousands times smaller, and the computation speed a few hundred times faster.

\begin{figure}[htbp]
\centering\includegraphics[width=1\textwidth]{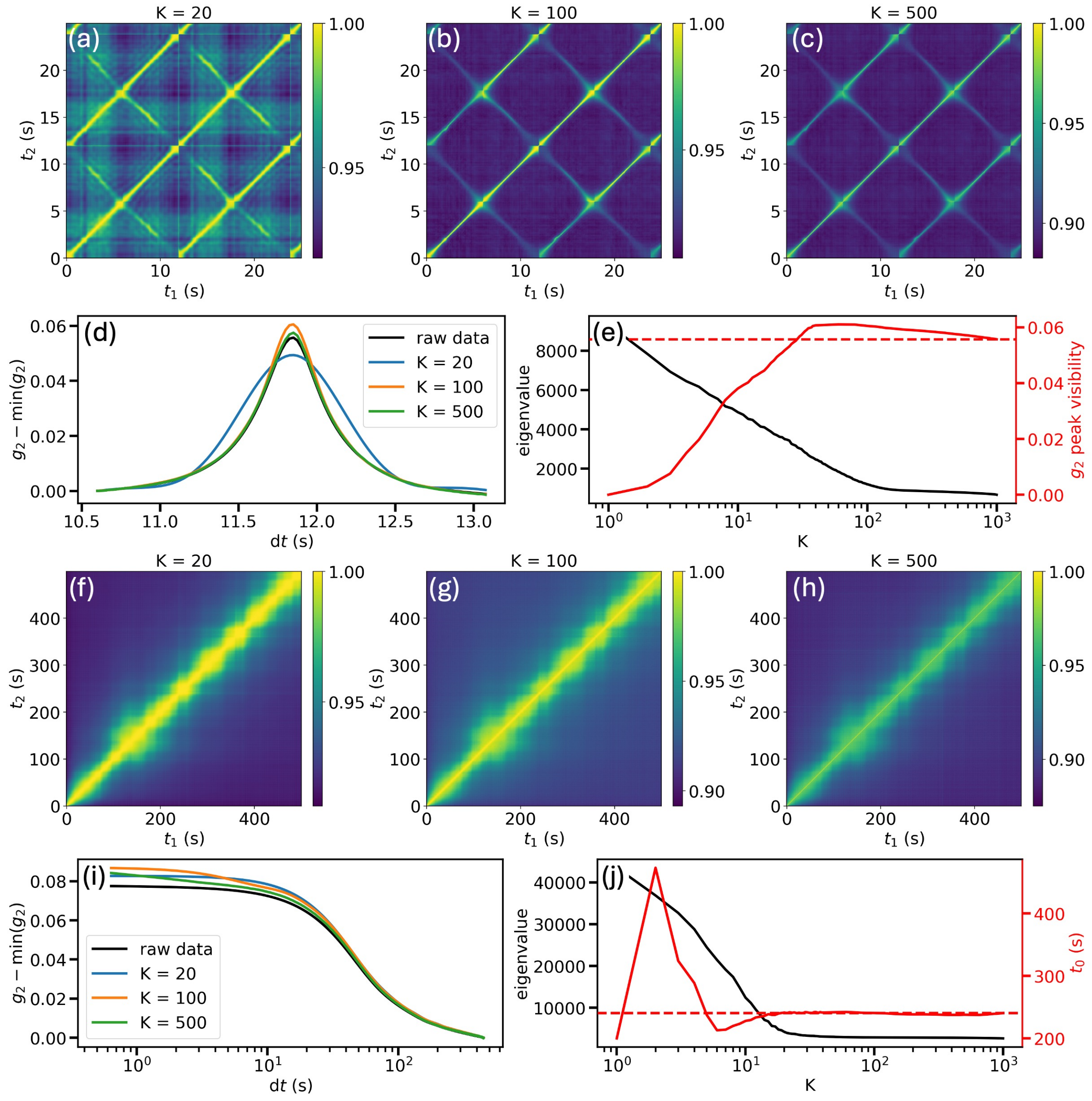}
\caption{Correlation calculation with offline lossy compression. TTC of lossy compressed oscillatory data with $K=$ (a) 20, (b) 100, and (c) 500. (d) shows a comparison of the peak at about $\text{d}t = 11.8$ s. (e) shows the eigenvalues of all the eigenvectors in \textbf{V} in descending order as well as the visibility of the peak shown in (d) as a function of $K$. TTC of lossy compressed rheology data with $K=$ (f) 20, (g) 100, and (h) 500. (i) shows a comparison of $g_2$. (j) shows  the eigenvalues of all the eigenvectors in \textbf{V} in descending order as well as the fitted relaxation time as a function of $K$.}
\label{fig:offline_lossy}
\end{figure}

Next we explore the scenario of lossy compression in the offline compression framework. With offline compression, the actual compression ratio is $\text{CR}\approx N/K$, which is the ratio between the size of the raw data $\textbf{X} \in \R^{N \times M}$ and $\textbf{V}_\text{K} \in \R^{M \times K}$. Fig.~\ref{fig:offline_lossy}a - c shows respectively the TTC for $K$ of 20, 100, and 500, corresponding to a CR of 50, 10, and 2. The oscillatory pattern can be clearly observed even at $K=20$. The use of lossy compression also results in a further reduction of the memory utilization and further acceleration of the computation, compared to the scenario of lossless compression. In Fig.~\ref{fig:offline_lossy}d we show a close-up look of the peak observed at about $\text{d}t = 11.8$ s in Fig.~\ref{fig:offline_lossless}c. We define peak visibility as the difference between the peak and the background level of the $g_2$ plot, and use it to quantify the fidelity of the TTC from the lossy compression. Fig.~\ref{fig:offline_lossy}e shows the peak visibility as a function of $K$. The peak visibility quickly improves with increasing $K$, approaching the level of the raw data (red dashed line) at about $K=30$. The peak visibility - $K$ relationship can be explained by the eigenvalues of the encoding matrix (black line). For this specific dataset, only 30 eigenvectors have eigenvalues that are twice as large as the minimum value. This indicates that the data itself can be reasonably described using a linear combination of just these 30 eigenvectors, and as such it is not surprising for the TTC of any $K>30$ compressed data to be a good approximation to the TTC of the raw data.

Similar conclusion can be drawn on the rheology data. Fig.~\ref{fig:offline_lossy}f - h shows respectively the TTC for $K$ of 20, 100, and 500, corresponding to a CR of 50, 10, and 2. The main features of the TTC are reproduced even at $K=20$. Fig.~\ref{fig:offline_lossy}i are the corresponding $g_2$ functions. The decorrelation of the $g_2$ function is reflective of the relaxation kinetics, and can be described using the Kohlrausch-Williams-Watts (KWW) expression\cite{Ross2010,Begam2021}.
\begin{equation}
g_2(\text{d}t) = B + C\;\text{exp}^{-2\text{d}t/t_0},
\label{eq:relaxation}
\end{equation}
where $B$ is the baseline, and $C$ is the contrast. $t_0$ is the relaxation time of the process. Fig.~\ref{fig:offline_lossy}j shows the fitted relaxation time as a function of $K$. The relaxation time of the raw data (red dashed line) is reproduced for $K>20$. This again can be understood using the eigenvalues of the encoding matrix (black line). For this specific dataset, only 20 eigenvectors have eigenvalues that are twice as large as the minimum value. This indicates that the data itself can be reasonably described using a linear combination of just these 20 eigenvectors, and as such it is not surprising for the $g_2$ of any $K>20$ compressed data to be good approximations to the $g_2$ of the raw data.

\subsection{Substantial Compression Ratio with Online Compression}\label{subsec32}
\begin{figure}[htbp]
\centering\includegraphics[width=1\textwidth]{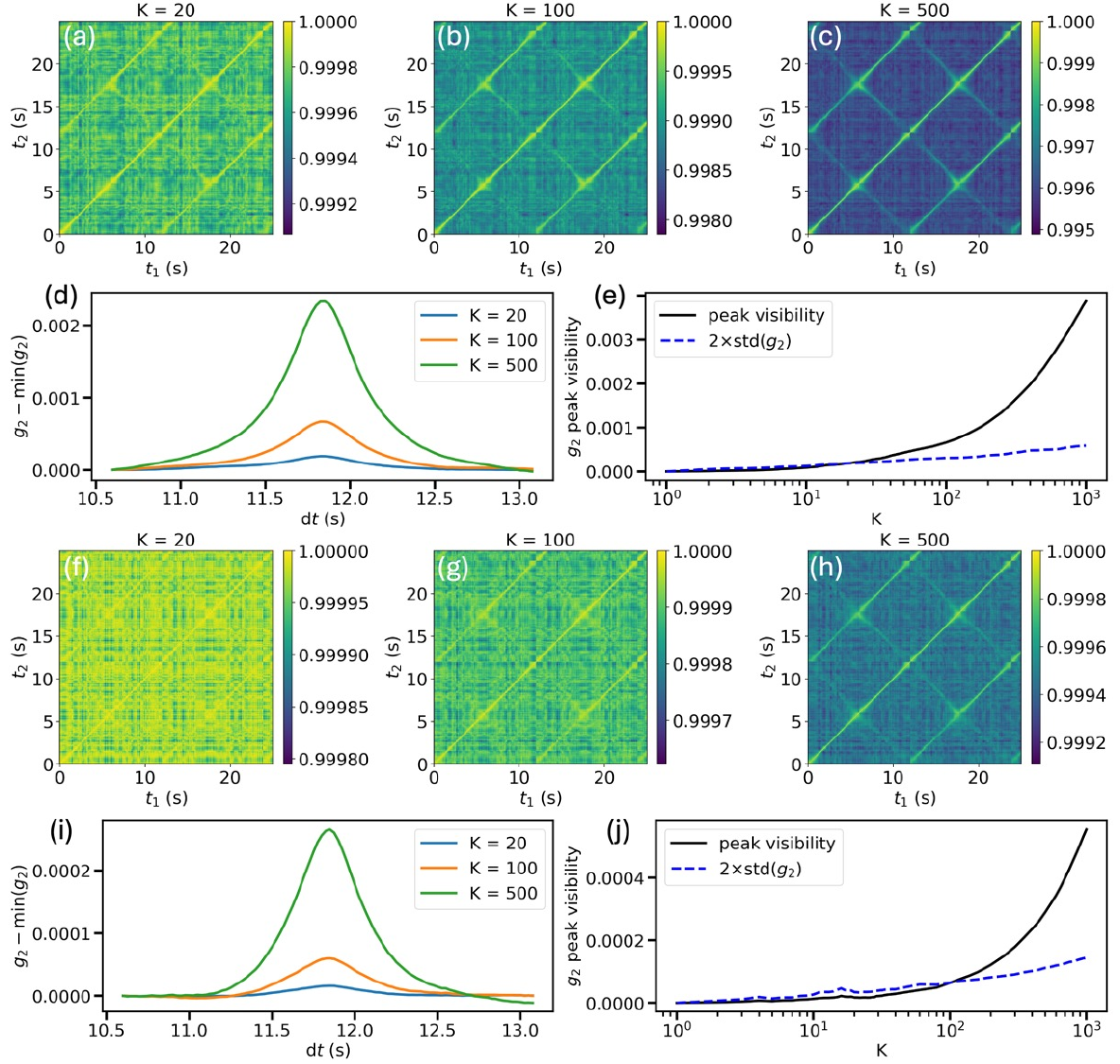}
\caption{Correlation calculation with online lossy compression. TTC of lossy compressed oscillatory data with $K=$ (a) 20, (b) 100, and (c) 500, using an encoding matrix generated on a related data. (d) shows a comparison of the peak at about $\text{d}t = 11.8$ s. (e) shows the visibility of the peak shown in (d) as well as the background level of the TTC as a function of $K$. TTC of lossy compressed oscillatory data with $K=$ (f) 20, (g) 100, and (h) 500, using an encoding matrix generated on an unrelated data. (i) shows a comparison of the peak at about $\text{d}t = 11.8$ s. (j) shows the visibility of the peak shown in (i) as well as the background level of the TTC as a function of $K$.}
\label{fig:online_lossy}
\end{figure}

For online compression, the encoding matrix $\textbf{V}_\text{K}$ is generated prior to the data acquisition and each image is compressed on the fly immediately after their individual acquisition. Because the same $\textbf{V}_\text{K}$ is used to compress all the data, the actual compression ratio becomes $\text{CR}\approx M/K$, which is the ratio between the size of the raw data $\textbf{X} \in \R^{N \times M}$ and the compressed data $\textbf{Y}_\text{K} \in \R^{N \times K}$. We also note that online compression is always lossy except for $K=\text{rank(\textbf{X})}$. Fig.~\ref{fig:online_lossy}a - c shows respectively the TTC for $K$ of 20, 100, and 500, corresponding to a substantial CR of 40000, 8000, and 1600. The encoding matrix in this case is generated using data acquired on the same sample during a previous measurement. Although the oscillatory pattern is still recognizable at $K=20$, the contrast is greatly reduced. Lossy SVD compression in general increases the similarity between the compressed data by removing small features and noise, which is why the computed correlation values in Fig.~\ref{fig:online_lossy}a - c are all close to 1. Similar effects are expected for lossy compression through quantization\cite{Huang2021} due to the limit number of possible data value. To quantify the loss of information with online compression, Fig.~\ref{fig:online_lossy}d shows the visibility of the peak for a $K$ value of 20, 100, and 500. Even for $K=500$, the absolute value of the peak visibility is about 30 times smaller than what is shown for the offline compression. Fig.~\ref{fig:online_lossy}e shows the peak visibility as a function of $K$. We consider the correlation peaks to be observable if it is stronger than twice the background level, calculated as the standard deviation of the TTC. With that, we find a minimum of $K=20$ is required in order for the oscillatory events to be detectable in this scenario.

It is also possible to generate the encoding matrix for the online compression using unrelated data. This approach has the advantage of $\textbf{V}_\text{K}$ being agnostic to the instrument and to the sample. In the following, the encoding matrix $\textbf{V}$ is generated using 1000 images after shifting the original mandrill image from the SIPI image database\cite{SIPI} randomly in the horizontal and vertical directions. For a given $K$ value, the encoding matrix $\textbf{V}_\text{K}$ for lossy compression uses the $K$ eigenvectors with the highest eigenvalues. Fig.~\ref{fig:online_lossy}f - h shows respectively the TTC for $K$ of 20, 100, and 500, corresponding to a CR of 40000, 8000, and 1600. The oscillatory pattern is barely recognizable at $K=20$. As shown in Fig.~\ref{fig:online_lossy}i, the peak visibility is about 10 times worse compared to what is shown in Fig.~\ref{fig:online_lossy}d. Using again twice the standard deviation as the cutoff, Fig.~\ref{fig:online_lossy}j shows that a larger $K$ of at least 100 is needed in this case for the oscillations to be considered detectable. Similar peak visibilities are obtained using other images ({e.g.}, the fishing boat) in the SIPI image database. We note that the idea of online compression is not to reproduce accurately the TTC or the $g_2$ function, but to aid in the detection of rare events in \insitu{}  experiments. The computation time for both the TTC and the $g_2$ increases quickly with the number of acquired images, and reaches tens of seconds for $N>10000$ images acquired on a megapixel detector ($M > 1000000$). This stands in contrast with the ever decreasing data acquisition time ({i.e.}, increasing frame rate) to the order of milliseconds, thanks to advent in detector technology and brighter light sources. Although the use of HPC resources can significantly reduce the computation time, on-demand resources are often scarce, and requires additionally a dedicated high-speed data transfer pipeline. With the online compression scheme proposed in this work, it is possible to compute the TTC and the $g_2$ at kHz framerate to provide feedback for rare event detection in real time. For the examples shown in this work, the encoding of a new image with $K=100$ takes 0.5 ms on an NVidia GeForce RTX 3090 GPU, and the subsequent computation of the TTC using the compressed data with $N=4000$ takes 1 ms. More importantly, the entire compressed data has a typical size of less than a few MB, which allows them to live on a small edge device for repeated real-time calculations. 


\section{Discussion}\label{sec4}

In this work, we show an approach towards direct operation on homomorphically compressed scientific data without their decompression. Two strategies have been proposed as illustrated in Fig.~\ref{fig:workflow}. The offline compression utilizes encoding matrix generated from the raw data itself. When the entirety of the encoding matrix is used, the compression is lossless. The photon correlation of the SVD compressed data is exactly the photon correlation of the raw data, despite the computation being performed on a dataset that is 800 times smaller and executes 800 times faster. To illustrate this, we have demonstrated on two sets of experimental data, one with an oscillatory correlation pattern, and the other with a relaxation kinetics. The second dataset is particularly intriguing, as it shows the possibility to compress and directly compute on a subsets of the 2D images for studying $q$ dependent behaviors\cite{Chen2016}. For (offline) lossy compression, the maximum compression ratio while retaining a good approximation of the raw data results can be inferred from the eigenvalues of the encoding matrix. A $K$ in the range of 20 to 30 is typically required for the representative data shown in this work, corresponding to an effective compression ratio of $\approx$40.

For online compression, a much larger compression ratio in the range of 1,000 to 10,000 can be achieved since the same encoding matrix is reused for the compression of all the data. When the encoding matrix is generated using data similar to the data to be compressed, qualitative approximation of the TTC and $g_2$ can be obtained with a relatively low $K$ value, because the eigenvectors in the encoding matrix captures some of the main features present in the raw data. It is also possible to generate the encoding matrix on completely unrelated data. The accuracy is about 10 times worse compared to the previous case, but the compression scheme becomes both instrument and sample agnostic. In either case, the effective computation time is greatly reduced to the order of milliseconds for each additional image. This enables real time calculation of the TTC and $g_2$ as new images are acquired at kHz framerate. More importantly, such calculation can be performed potentially on an edge device to providing timely feedback for experimental steering purposes, without the necessity of on-demand HPC resources.

Our approach offers an effective solution to the big data challenge stemming from the drastically increased raw data rate at coherent light sources. While there have been extensive research on how to effectively compress scientific data, there has been no demonstration, to our knowledge, that generates quantitatively accurate results directly using compressed dataset. The data processing for XPCS is in particular challenging because computation needs to be performed on the entire dataset (and not just the latest image) that grows over time. The complexity for computing the photon correlation using direct matrix multiplication is $\bigO(n^2m)$, where $n=N$ is the number of images in the time series, $m = M$ for the raw data and $m=K$ for the compressed data. This explains why a $K=100$ can lead to substantial acceleration of the compute for data acquired on a megapixel detector ($M>1,000,000$). 

Beyond the example demonstrated in this work, our framework can be extended to facilitate real time operation directly on a compressed data stream for other techniques. The ideal compression algorithms will be different for different techniques as they need to be tailored to the corresponding data analysis routines, but it is recommended that they share the following characteristics as the SVD compression used in this work. First, a fixed-length compressed data is produced which simplifies any matrix operation that includes both computer vision and neural network based approaches. Second, the lossy compression of the dataset is achieved while retaining the most significant variance and reducing noise. More importantly, the compression highlights the most influential components in the data, thus facilitating a more interpretable understanding of the underlying patterns and relationships.




\backmatter


\section*{Acknowledgements}
Work performed at the Center for Nanoscale Materials and Advanced Photon Source, both U.S. Department of Energy Office of Science User Facilities, was supported by the U.S. DOE, Office of Basic Energy Sciences, under Contract No. DE-AC02-06CH11357. This work was supported by the U.S. Department of Energy, Office of Science, Advanced Scientific Computing Research and the Basic Energy Sciences Accelerator and Detector R\&D program, under Contract DE-AC02-06CH11357. Part of the inception of ideas and analysis was supported by the U.S. Department of Energy, Office of Science, Basic Energy Sciences, Materials Science and Engineering Division through the Early Career Research Program. We thank Dr. Panpan Huang and Prof. Chris Jacobsen for the inspiring conversations related to data compression. We thank Mike Hammer, Henry Shi, Senthil Gnanasekaran and John Weizeorick for stimulating discussions related to hardware co-design. We acknowledge Prof. Mark Foster, University of Akron and Prof. Robert Leheny, Johns Hopkins for giving permission to use their data.

\section*{Declarations}

\begin{itemize}
\item Funding

U.S. DOE, Office of Basic Energy Sciences Contract No. DE-AC02-06CH11357.

\item Conflict of interest/Competing interests (check journal-specific guidelines for which heading to use)

The authors declare no competing interests.

\item Ethics approval and consent to participate

Not Applicable.

\item Consent for publication

Authors give consent for publication.

\item Data availability 

The script for the analysis used in this work will be available on a public GitHub page.

\item Materials availability

Not Applicable.

\item Code availability 

The script for the analysis used in this work will be available on a public GitHub page.

\item Author contribution

S.S., A.M. and T.Z. conceived the study. S.S., Z.W.D., K.Y. and T.Z. contributed to the mathematical part of the project. S.S., M.C., M.V.H and T.Z. contributed to the computational part of the project. Y.C., Q.Z., E.M.D. and S.N. contributed to the photon science part of the project. S.S. and T.Z. performed the data analysis and wrote the manuscript, with input from all authors.
\end{itemize}


\end{document}